\documentclass{amsart}
\usepackage{graphicx}
\usepackage{amssymb}
\usepackage{amsfonts}
\swapnumbers
\newtheorem{theorem}{Theorem}[section]

\newtheorem{corollary}[theorem]{Corollary}

\theoremstyle{definition}

\newtheorem{assumption}[theorem]{Assumption}

\numberwithin{equation}{section}
 \theoremstyle{plain}
 
 \numberwithin{equation}{section} 
 \numberwithin{figure}{section} 
 \theoremstyle{plain}
 \theoremstyle{remark}
 \newtheorem*{acknowledgement*}{Acknowledgement}

\newcommand{\cA}{{\mathcal A}}

\newcommand{\cC}{{\mathcal C}}

\newcommand{\cF}{{\mathcal F}}
\newcommand{\cG}{{\mathcal G}}
\newcommand{\cH}{{\mathcal H}}

\newcommand{\cL}{{\mathcal L}}

\newcommand{\cN}{{\mathcal N}}

\newcommand{\Om}{{\Omega}}
\newcommand{\om}{{\omega}}
\newcommand{\ve}{{\varepsilon}}
\newcommand{\del}{{\delta}}
\newcommand{\Del}{{\Delta}}
\newcommand{\gam}{{\gamma}}
\newcommand{\Gam}{{\Gamma}}

\newcommand{\vr}{{\varrho}}
\newcommand{\Sig}{{\Sigma}}
\newcommand{\sig}{{\sigma}}
\newcommand{\al}{{\alpha}}
\newcommand{\be}{{\beta}}
\newcommand{\ka}{{\kappa}}
\newcommand{\la}{{\lambda}}

\newcommand{\up}{{\upsilon}}
\newcommand{\vrho}{{\varrho}}

\newcommand{\bbI}{{\mathbb I}}

\newcommand{\bfY}{{\bf Y}}
\newcommand{\bfX}{{\bf X}}

\begin{document}
\title[]{Limit theorems for numbers of returns\\
in arrays under $\phi$-mixing}%
 \vskip 0.1cm
 \author{ Yuri Kifer\\
\vskip 0.1cm
 Institute  of Mathematics\\
Hebrew University\\
Jerusalem, Israel}%
\address{
Institute of Mathematics, The Hebrew University, Jerusalem 91904, Israel}
\email{ kifer@math.huji.ac.il}%

\thanks{ }
\dedicatory{}
 \date{\today}
\begin{abstract}\noindent
We consider a $\phi$-mixing shift $T$ on a sequence space $\Om$ and study the number $\cN_N$ of returns $\{ T^{q_N(n)}\om\in A^a_n\}$
at times $q_N(n)$ to a cylinder $A^a_n$ constructed by a sequence $a\in\Om$ where $n$ runs either until a fixed integer $N$ or
 until a time $\tau_N$ of the first return $\{ T^{q_N(n)}\om\in A^b_m\}$ to another cylinder  $A^b_m$ constructed by $b\in\Om$. Here $q_N(n)$ are
 certain functions of $n$ taking on nonnegative integer values when $n$ runs from 0 to $N$ and the dependence on $N$ is the main generalization
 here in comparison to \cite{KY}. Still, the dependence on $N$ requires certain conditions under which we obtain Poisson distributions
 limits of $\cN_N$ when counting is until $N$ as $N\to\infty$ and geometric distributions limits when counting is until $\tau_N$ as $N\to\infty$.
 The results and the setup are similar to \cite{Ki19} where multiple returns are considered but under the stronger $\psi$-mixing assumption.
\end{abstract}
\maketitle
\markboth{Yu.Kifer}{Numbers of returns}
\renewcommand{\theequation}{\arabic{section}.\arabic{equation}}
\pagenumbering{arabic}

\section{Introduction}\label{sec1}\setcounter{equation}{0}

The study of returns to (hits of) shrinking targets by a dynamical system,
started in \cite{Pi}, \cite{Hi} and \cite{De}, has already about 30 years
history. These works were extended in various directions, in particular, to
returns to shrinking geometric balls by uniformly and non-uniformly hyperbolic
dynamical systems (see, for instance, \cite{HW} and references there), to multiple
returns to shrinking cylinder sets under $\psi$-mixing (see, for instance, \cite{KR1})
and others. More recently, motivated by the research on open dynamical systems (see, for instance,
\cite{DWY}) the
asymptotic behaviour of numbers of returns to a shrinking target until the first arrival
to another shrinking target was investigated in \cite{KR2} and \cite{KY} where the first
work dealt with the $\psi$-mixing case while the second one dealt with a $\phi$-mixing
situation which allowed applications to a wider class of dynamical systems. Another
generalization started in Ch.3 of \cite{HK} and continued in \cite{Ki19} and \cite{Ha},
dealt with returns at prescribed times which depended also on the total observation time where
additional peculiarities appeared.

In this paper we consider two related types of limit theorems for numbers of returns
which are represented by the sums
\[
S_N=\sum_{k=1}^N\bbI_{A^a_{n_N}}\circ T^{q_N(k)}\quad\mbox{and}\quad\Sig_N=\sum_{k=1}^{\tau_N}\bbI_{A^a_{n_N}}\circ T^{q_N(k)}
\]
where $\bbI_\Gam$ is the indicator of a set $\Gam$, $A^a_n$ is a cylinder set of the length $n$ constructed by a sequence $a$,
$T$ is a $\phi$-mixing left shift on a sequence space $\Om$, $\tau_N(\om)$ is the first $k$ such that $T^{q_N(k)}\om$ belongs
to another cylinder set $A^b_{m_N}$ and $q_N(k)$ are certain functions taking on nonnegative integer values when $k$ runs from 1
to $N$ for $N=1,2,...$. In probability such sums where summands themselves depend on the number of summands are usually called
(triangular) arrays. We will provide conditions on functions $q_N(k)$ such that as $N\to\infty$ the sum $S_N$ converges in
distribution to a Poisson random variable while the sum $\Sig_N$ converges in distribution to a geometric random variable. It is
easy to see that without certain conditions such results do not hold, in general. Indeed, taking $q_N(k)=k(N-k)$ we obtain that
the above sums may converge only to a random variable taking on just even values, and so the limits cannot have Poisson or geometric
distributions.

Our results remain valid for dynamical systems possessing appropriate symbolic representations such as Axiom A
 diffeomorphisms (see \cite{Bo}), expanding transformations and some maps having symbolic representations with an
 infinite alphabet and a $\psi$-mixing invariant measure such as the Gauss map with its Gauss invariant measure and
 more general $f$-expansions (see \cite{He}). A direct application of the above results in the symbolic setup yields
 the corresponding results for arrivals to elements of Markov partitions but employing additional technique (see, for
 instance, \cite{KY}) it is not difficult to extend these results for arrivals to shrinking geometric balls.
 Since we assume only $\phi$-mixing, rather than $\psi$-mixing, our results remain valid for some classes of
 nonuniformly expanding maps of the interval such as Gibbs-Markov maps and some others (cf. \cite{KY}). In the probability
 direction we can consider Markov chains with countable state spaces satisfying the Doeblin condition which are known to be
 exponentially fast $\phi$-mixing (see \cite{Br}), and so our results are applicable to the corresponding shifts in the path
 spaces.

 \section{Preliminaries and main results}\label{sec2}\setcounter{equation}{0}

Our setup consists of a finite or countable set $\mathcal{A}$ which is not a singleton,
the sequence space $\Omega=\mathcal{A}^{\mathbb{N}}$,
the $\sigma$-algebra $\mathcal{F}$ on $\Omega$ generated by cylinder
sets, the left shift $T:\Omega\rightarrow\Omega$, and a $T$-invariant
probability measure $P$ on $(\Omega,\mathcal{F})$ which is assumed to be
$\phi$-mixing with respect to the $\sig$-algebras $\cF_{mn},\, n\geq m$ generated
by the cylinder sets of the form $\{\om=(\om_0,\om_1,...)\in\Om:\,\om_i=a_i\,$ for
$m\leq i\leq n\}$ for some $a_m,a_{m+1},...,a_n\in\cA$. Observe also that
$\cF_{mn}=T^{-m}\cF_{0,n-m}$ for $n\geq m$.

Recall, that the $\phi$-dependence (mixing) coefficient between two $\sig$-algebras
$\cG$ and $\cH$ can be written in the form (see \cite{Br}),
\begin{eqnarray}\label{2.1}
&\phi(\cG,\cH)=\sup_{\Gam\in\cG,\Del\in\cH}\big\{\big\vert\frac {P(\Gam\cap\Del)}
{P(\Gam)}-P(\Del)\big\vert,\, P(\Gam)\ne 0\big\}\\
&=\frac 12\sup\{\| E(g|\cG)-E(g)\|_{L^\infty}:\, g\,\,\mbox{is}\,\,
\cH-\mbox{measurable and}\,\, \|g\|_{L^\infty}\leq 1\}.\nonumber
\end{eqnarray}
Set also
\begin{equation*}
\phi(n)=\sup_{m\geq 0}\phi(\cF_{0,m},\cF_{m+n,\infty}).
\end{equation*}
The probability $P$ is called $\phi$-mixing if $\phi(n)\to 0$ as $n\to\infty$.

We will need also the $\al$-dependence (mixing) coefficient between two $\sig$-algebras
$\cG$ and $\cH$ which can be written in the form (see \cite{Br}),
\begin{eqnarray}\label{2.2}
&\al(\cG,\cH)=\sup_{\Gam\in\cG,\Del\in\cH}\big\{\big\vert P(\Gam\cap\Del)
-P(\Gam)P(\Del)\big\vert\big\}\\
&=\frac 14\sup\{\| E(g|\cG)-E(g)\|_{L^1}:\, g\,\,\mbox{is}\,\,
\cH-\mbox{measurable and}\,\, \|g\|_{L^\infty}\leq 1\}.\nonumber
\end{eqnarray}
Set also
\begin{equation*}
\al(n)=\sup_{m\geq 0}\al(\cF_{0,m},\cF_{m+n,\infty}).
\end{equation*}

 For each word $a=(a_0,a_1,...,a_{n-1})\in
\cA^n$ we will use the notation $[a]=\{\om=(\om_0,\om_1,...):\, \om_i=a_i,\,
 i=0,1,...,n-1\}$ for the corresponding cylinder set. Without loss of generality
 we assume that the probability of each 1-cylinder set is positive, i.e. $P([a])>0$
 for every $a\in\cA$, and since $\cA$ is not a singleton we have also
 $\sup_{a\in\cA}P([a])<1$. Write $\Omega_{P}$ for the support
of $P$, i.e.
\[
\Om_P=\{\om\in\Om:\,P[\om_0,...,\om_n]>0\,\,\mbox{for all}\,\, n\geq 0\}.
\]
 For $n\ge 1$ set $\mathcal{C}_{n}=\{[w]\::\:w\in\mathcal{A}^{n}\}$. Then $\cF_{0,n}$
 consists of $\emptyset$ and all unions of disjoint elements from $\cC_{n+1}$. If the probability
 $P$ is $\phi$-mixing then by Lemma 3.1 from \cite{KY} there exists $\up>0$ such that
\begin{equation}\label{2.3}
P(A)\le e^{-\upsilon n}\text{ for all \ensuremath{n\ge1} and
\ensuremath{A\in\mathcal{C}_{n}}}.
\end{equation}

Next, for any $U\in\cF_{0,n-1},\, U\ne\emptyset$ define
\[
\pi(U)=\min\{k\geq 1:\:U\cap T^{-k}U\ne\emptyset\}
\]
and observe that $\pi(U)\leq n$. We will be counting the returns to $U$ at times $q_N(k)$
considering the sum
\[
S^U_N=\sum_{k=1}^N\bbI_U\circ T^{q_N(k)}.
\]
Our goal will be to show that if $U$ is replaced by a sequence of Borel sets $U_N\subset\Om$
such that $NP(U_N)$ converges as $N\to\infty$ then $S_N^{U_N}$ converges in distribution to a
Poisson random variable and, as an example in Introduction shows, in order to achieve this some
assumptions on functions $q_N(k)$ are necessary.
\begin{assumption}\label{ass2.1} $q_N(n)$ is a function taking on nonnegative integer values
 on integers $n,N\geq 0$, defined arbitrarily when $n>N$ and such that for some constant $K>0$ and all $N\geq 1$ the
 following properties hold true:

 (i) For all $k$ the number of integers $n,\, 0\leq n\leq N$ satisfying the equation
  \begin{equation*}
 q_N(n)=k
 \end{equation*}
 does not exceed $K$;

 (ii) The number of pairs $m\ne n$ satisfying $0\leq m,n\leq N$ and solving the equation
 \begin{equation*}
 q_N(n)-q_N(m)=0
 \end{equation*}
 does not exceed $K$;
\end{assumption}

First, note that the example of $q_N(k)=k(N-k)$ from Introduction does not satisfy Assumption \ref{ass2.1}(ii)
since $q_N(n)=q_N(N-n)$, and so at least $[N/2]-1$ pairs $n\ne m=N-n$ solve the equation in (ii).
Next, observe that if there exist $n_0,N_0\geq 1$ such that for all $N\geq N_0$ the function $q_N(n)$ of $n$ is
strictly increasing when $n_0\leq n\leq N$, then the whole Assumption \ref{ass2.1} is satisfied. Indeed, at most
one $n\geq n_0$ can solve the equation $q_N(n)=k$ when $N$ and $k$ are fixed, and so the number of solutions in (i)
cannot exceed $n_0+1$. Next, if $N\geq n,m\geq n_0$ then $q_N(n)=q_N(m)$ will hold true only if $n=m$. If, say, $m<n_0$
then for such $m$ there could be at most one $n\geq n_0$ satisfying $q_N(n)=q_N(m)$. It follows that there exist at most
$n_0^2+2n_0$ pairs $0\leq m,n\leq N$ such that $q_N(n)=q_N(m)$.
 In particular, if $q_N(n)=r(n)+g(N)$ where $r$ is a nonconstant polynomial in $n$
 and $g$ is a function of $N$, both nonnegative for $n,N\geq 0$ and taking on integer values on integers, then $q_N$
 satisfies Assumption \ref{ass2.1}. Indeed, the number of solutions in Assumption \ref{ass2.1}(i) is bounded by the degree
 of $r$ and there exists an integer $n_0\geq 1$ such that the polynomial $r$ is strictly increasing on $[n_0,\infty)$.

  For any two random variables or random vectors $Y$ and $Z$ of the same
dimension denote by $\cL(Y)$ and $\cL(Z)$ their distribution and by
\[
d_{TV}(\cL(Y),\,\cL(Z))=\sup_G|\cL(Y)(G)-\cL(Z)(G)|
\]
the total variation distance between $\cL(Y)$ and $\cL(Z)$ where the supremum
is taken over all Borel sets. Denote by Pois$(\la)$ the Poisson distribution with
a parameter $\la>0$, i.e. Pois$(\la)(k)=e^{-\la}\frac {\la^k}{k!}$ for each $k=0,1,2,...$.
Our first result is the following.

\begin{theorem}\label{thm2.2} Suppose that Assumption \ref{ass2.1} is satisfied.  Then there exists a
constant $C\geq 1$ such that for any $n,\, V\in\cF_{0,n-1}$, $N$ and $R$,
 \begin{eqnarray}\label{2.4}
 &d_{TV}(\cL(S_N^V),\,\mbox{Pois}(\la_N))\leq CN\big(R(P(V))^2\\
 &+P(V)\sum_{r=\pi(V)}^R(\phi([r/2]+1)+P(T^{n-[r/2]}V))+\phi(R-n)\big)\nonumber
 \end{eqnarray}
 where $\la_N=N(P(V))^\ell$.
\end{theorem}
We observe that related estimates under $\phi$-mixing were obtained in \cite{AV} but it is difficult to obtain definitive
convergence results from there.
 \begin{corollary}\label{cor2.3} Suppose that Assumption \ref{ass2.1} is satisfied and the $\phi$-mixing coefficient
 is summable, i.e.
 \[
 \sum_{k=1}^\infty\phi(k)<\infty.
 \]
 Let $V_{L}\in\cF_{0,n_L-1},\, L=1,2,...$ be a sequence of sets such that
 $n_LP(V_{L})\to 0$ and $\sum_{r=\pi(V_{L})}^{n_L-1}P(T^{n_L-r}V_{L})\to 0$ as $L\to\infty$. Let $N_L\to\infty$
 as $L\to\infty$ be a sequence of integers such that $0<C^{-1}\leq\la_{L}=N_LP(V_{L})\leq C<\infty$ for some constant
  $C$ and all $L\geq 1$. Then
 \begin{equation}\label{2.5}
 d_{TV}(\cL(S_{N_L}^{V_{L}}),\,\mbox{Pois}(\la_{L}))\to 0\,\,\mbox{as}\,\, L\to\infty
 \end{equation}
 and if $\lim_{L\to\infty}\la_{L}=\la$ then the distribution of $S_{N_L}^{V_{N_L}}$ converges in total
 variation as $L\to\infty$ to the Poisson distribution with the parameter $\la$. In particular, if
 $V_{L}=A^\eta_{n_L}=[\eta_0,...,\eta_{n_L-1}]=\{\om\in\Om:\,\om_0=\eta_0,...,\om_{n_L-1}=\eta_{n_L-1}\}$ with $n_L\to\infty$
 as $L\to\infty$ and $\eta\in\Om_P$ is nonperiodic then $\pi(A^\eta_{n_L})\to\infty$ as $L\to\infty$ and the above
 statements hold true for such $V_{L}$'s provided the above conditions on $\la_L$ are satisfied.
 \end{corollary}

 Next, for any $V\in\cF_{0,n-1},\, V\ne\emptyset$ and $W\in\cF_{0,m-1},\, W\ne\emptyset$ define
\[
\pi(V,W)=\min\{k\geq 1:\:V\cap T^{-k}W\ne\emptyset
\mbox{ or }W\cap T^{-k}V\ne\emptyset\}\:.
\]
It is clear that $\pi(V,W)\leq m\wedge n$, and so
\[
\kappa_{V,W}=\min\{\pi(V,W),\pi(V),\pi(W)\}\leq m\wedge n
\]
where, as usual, for $n,m\ge1$ we denote $m\vee n=\max\{m,n\}$ and $m\wedge n=\min\{m,n\}$.
Set
\[
\tau_W(\omega)=\min\{k\ge 1\::\:T^{q_N(k)}\omega\in W\}
\]
with $\tau_W(\om)=\infty$ if the event in braces does not occur and define
\[
\Sig_N^{V,W}=\sum_{k=1}^{\tau_W}\bbI_{V}\circ T^{q_N(k)}.
\]
Denote by Geo$(\rho)$ the geometric distribution with a parameter $\rho\in(0,1)$,
i.e. Geo$(\rho)(k)=\rho(1-\rho)^k$ for each $k=0,1,2,...$.

\begin{theorem}\label{thm2.4} Assume that Assumption \ref{ass2.1} is
satisfied. Then there exists a constant $C>0$ such that for any disjoint sets
 $V\in\cF_{0,n-1}$ and $W\in\cF_{0,m-1}$ with $P(V),P(W)>0$ and all integers $n,m,N,R\ge 1$,
\begin{eqnarray}\label{2.6}
&\quad\quad d_{TV}(\mathcal{L}\big(\Sig_N^{V,W}),Geo(\rho)\big)\le C\bigg((1-P(W))^N+(n\vee m)(P(V)+P(W))\\
&+RN(P(V)+P(W))^2+N\phi(R-n\vee m)\nonumber\\
&+N(P(V)+P(W))\sum_{r=\ka_{V,W}}^{n\vee m-1}(\phi(r)+P(T^{n\vee m-r}V)+(P(T^{n\vee m-r}W)\big)\bigg)
\nonumber\end{eqnarray}
where $\rho=\frac{P(W)}{P(V)+P(W)}$.
\end{theorem}

\begin{corollary}\label{cor2.5}
Suppose that Assumption \ref{ass2.1} holds true and the $\phi$-mixing coefficient
is summable. Let $V_L\in\cF_{0,n_L-1}$ and $W_L\in\cF_{0,m_L-1}$, $L=1,2,...$
be two sequences of sets such that
\begin{equation}\label{2.7}
(n_L\vee m_L)(P(V_L)+P(W_L))\to 0,\,\,\ka_{V_L,W_L}\to\infty\quad\mbox{as}\quad L\to\infty,
\end{equation}
\begin{equation}\label{2.8}
\al_L=\sum_{r=\ka_{V_L,W_L}}^{n_L\vee m_L-1}(P(T^{n_L\vee m_L-r}V_L)+P(T^{n_L\vee m_L-r}W_L))\to 0\,\,\mbox{as}\,\, L\to\infty
\end{equation}
and for some constant $C$ and all $L\geq 1$,
\begin{equation}\label{2.9}
 0<C^{-1}\leq\frac {P(V_L)}{P(W_L)}\leq C<\infty.
\end{equation}
Let $N_L,\, L=1,2,...$ be a sequence satisfying
\begin{equation}\label{2.10}
N_LP(W_L)\to\infty\,\,\mbox{and}\,\, N_L(M_{N_L}+n_L\vee m_L+\al_L)(P(W_L))^2\to 0\,\,\mbox{as}\,\, L\to\infty
\end{equation}
where $M_N=M_N^{(\ve)}=\min\{ n\geq 1:\,\frac n{\gam^\ve(n)}\geq N\}$ for some $0<\ve <1$ and $\gam(n)=n\phi(n)$.
Then
\begin{equation}\label{2.12}
d_{TV}(\cL(\Sig_{N_L}^{V_L,W_L}),\, Geo(\rho_L))\to 0\,\,\mbox{as}\,\, L\to\infty
\end{equation}
where $\rho_L=P(W_L)(P(W_L)+P(V_L))^{-1}$. In particular, if $\lim_{L\to\infty}\rho_L=\rho$,
then $\Sig_{N_L}^{V_L,W_L}$ converges in total variation as $L\to\infty$ to the geometric distribution with the parameter
$\rho$. Furthermore, let $V_L=A^\xi_{n_L}=[\xi_0,...,\xi_{n_L-1}]\in\cC_{n_L}$ and $W_L=A^\eta_{m_L}=[\eta_0,...,\eta_{m_L-1}]\in\cC_{m_L}$
with $n_L,m_L\to\infty$ as $L\to\infty$ and suppose that $\xi,\eta$ are not periodic and not shifts of each other. Then
\begin{equation}\label{2.13}
\ka_{A^\xi_{n_L},A^\eta_{m_L}}\to\infty\,\,\mbox{as}\,\, L\to\infty
\end{equation}
and if also
\begin{equation}\label{2.14}
n_L\wedge m_L+\ka_{A^\xi_{n_L},A^\eta_{m_L}}-n_L\vee m_L\to\infty\,\,\mbox{as}\,\, L\to\infty
\end{equation}
then (\ref{2.8}) holds true. In fact, (\ref{2.14}) is satisfied for $P\times P$-almost all $(\xi,\eta)\in\Om\times\Om$
provided
\begin{equation}\label{2.15}
2n_L\wedge m_L-n_L\vee m_L-3\up\ln(n_L\wedge m_L)\to\infty\,\,\mbox{as}\,\, L\to\infty
\end{equation}
where $\up$ is from (\ref{2.3}).
\end{corollary}

Observe that when $q_N(n)$ does not depend on $N$ then $\Sig_N^{V,W}$ does not depend on $N$ either and in order to
obtain (\ref{2.12}) relying on (\ref{2.6}) we have only to pick up some sequence $N_L$ satisfying (\ref{2.10}) which is
always possible provided (\ref{2.7})--(\ref{2.9}) hold true.

\section{Poisson distribution limits}\label{sec3}\setcounter{equation}{0}

We will need the following semi-metrics between positive integers $k,l> 0$,
\[
\del_N(k,l)=|q_{N}(k)-q_{N}(l)|.
\]
It follows from Assumption \ref{ass2.1}(i) that for any integers $n\in\{ 1,...,N\}$ and $k\geq 0$,
\begin{equation}\label{3.3}
\#\{ m:\, \del_N(n,m)=k\}\leq 2K.
\end{equation}
For any integers $M,R\geq 1$ and $1\leq n\leq N$ introduce the sets
\[
B^{M,R}_{n,N}=\{ l:\, 1\leq l\leq M,\,\del_N(l,n)<R\}\quad\mbox{and}\quad B^R_{n,N}=B^{N,R}_{n,N}.
\]
By (\ref{3.3}), for any $n$,
\begin{equation}\label{3.4}
\# B^{M,R}_{n,N}\leq\min(M,\, 2KR).
\end{equation}

Let $V\in\cF_{0,n-1}$ and set $X_{k,N}=X^V_{k,N}=\bbI_V\circ T^{q_{N}(k)}$.
Then $S_N=S_N^V=\sum_{k=1}^NX_{k,N}$. Set $p_{k,N}=P\{ X_{k,N}=1\}$ and $p_{k,l,N}=P\{ X_{k,N}=1$ and $X_{l,N}=1\}$.
Since $T$ is $P$-preserving $p_{k,N}=E(\bbI_V\circ T^{q_{N}(k)})=P(V)$ and $p_{k,l,N}=P(V\cap T^{-(q_{N}(k)-q_N(l))}V)$
provided $q_N(l)\leq q_N(k)$. By Theorem 1 from \cite{AGG} we obtain
\begin{equation}\label{6.1}
d_{TV}(\cL(S_N),\, Pois(\la_N))\leq b_1+b_2+b_3
\end{equation}
where $b_1,\, b_2$ and $b_3$ are defined by
\begin{equation}\label{4.2}
b_1=\sum_{n=1}^N\sum_{l\in B^R_{n,N}}p_{n,N}p_{l,N},\,\,\, b_2=\sum_{n=1}^N\sum_{n\ne l\in B^R_{n,N}}p_{n,l,N}
\end{equation}
and
\begin{equation}\label{4.3}
b_3=\sum_{n=1}^Ns_{n,N}\,\,\mbox{with}\,\, s_{n,N}=E|E(X_{n,N}-p_{n,N}|\sig\{ X_{l,N}:\, l\in\{ 1,...,N\}\setminus
B^R_{n,N}\})|.
\end{equation}
By (\ref{3.4}) and (\ref{4.2}) we conclude that
\begin{equation}\label{6.4}
b_1=\sum_{k=1}^N\sum_{l\in B^R_{k,N}}p_{k,N}p_{l,N}\leq 2KRN(P(V))^2.
\end{equation}

In order to estimate $p_{k,l,N}$ we observe that if $|i-j|<\pi(V)$ then $(\bbI_V\circ T^i)(\bbI_V\circ T^j)=0$. Hence,
$p_{k,l,N}=0$ if $\del_N(k,l)<\pi(V)$.
Now suppose that $\del_N(k,l)=d$ with $\pi(V)\leq d<n$. Then
\begin{equation}\label{6.5}
\mbox{either}\,\,\, q_{N}(l)\leq q_{N}(k)-d\,\,\,\mbox{or}\,\,\, q_{N}(l)\geq q_{N}(k)+d.
\end{equation}
 Assume, for instance, that the first inequality in (\ref{6.5}) holds true and let $r=q_{N}(k)-q_{N}(l)$.
Then $r\geq d\geq\pi(V)$. If $r\geq n$ then since $V\in\cF_{0,n-1}$, we obtain by the definition of the $\phi$-mixing
coefficient that
\begin{equation}\label{6.6}
p_{k,l,N}=P(V\cap T^{-r}V)\leq(\phi(r-n+1)+P(V))P(V).
\end{equation}
Suppose that $\pi(V)\leq r<n$ and assume that $V\cap T^{-r}V\ne\emptyset$. Let $s\geq n-r$, set $V_s=T^sV$ and observe that
$T^{-s}V_s\supset V$. Then by the definition of the $\phi$-mixing coefficient,
\begin{eqnarray}\label{6.7}
&p_{k,l,N}=P(V\cap T^{-r}V)\leq P(V\cap T^{-(r+s)}V_s)\leq (\phi(r+s-n+1)\\
&+P(T^sV))P(V)\leq (\phi([r/2]+1)+P(T^{n-[r/2]}V))P(V)\nonumber
\end{eqnarray}
taking $s= n-[r/2]$. If the second inequality in (\ref{6.5}) holds true then we obtain (\ref{6.6}) if $r=q_{N}(l)-q_{N}(k)\geq n$,
  while if $\pi(V)\leq r<n$ then we arrive at (\ref{6.7}). Observe that by Assumption \ref{ass2.1}(i) for any $N\geq 1$
  and integers $k\geq 0$ and $r$,
  \begin{equation}\label{6.9}
  \#\{ l\geq 0:\, q_{N}(k)-q_{N}(l)=r\}\leq K.
  \end{equation}
  Now, it follows from (\ref{3.3}), (\ref{3.4}) and (\ref{6.6})--(\ref{6.9}) that
  \begin{equation}\label{6.10}
b_2=\sum_{k=1}^N\sum_{k\ne l\in B_{k,N}^R}p_{k,l,N}\leq 4KNP(V)\sum_{r=\pi(V)}^R(\phi([r/2]+1)+P(T^{n-[r/2]}V)).
\end{equation}

Next, we estimate $s_{k,N}$ and $b_3$ defined by (\ref{4.3}). Since $\del_N(k,l)\geq R$ for $l\not\in B^R_{k,N}$ and $V\in\cF_{0,n-1}$,
we derive from Lemma 3.3 in \cite{KY} and the definition of the $\al$-mixing coefficient that for $n<R<N$,
\begin{equation}\label{6.11}
s_{k,N}\leq\al\big(\cF_{q_N(k),q_N(k)+n},\sig(\cF_{0,q_N(k)-R+n},\cF_{q_N(k)+R-n,\infty})\big)\leq 3\phi(R-n).
\end{equation}
Hence, by (\ref{4.3}) and (\ref{6.11}),
\begin{equation}\label{6.13}
b_3=\sum_{k=1}^Ns_{k,N}\leq 3N\phi(R-n).
\end{equation}
Finally, collecting (\ref{6.1}), (\ref{6.4}), (\ref{6.10}) and (\ref{6.13}) we derive (\ref{2.4}) completing the proof of
Theorem \ref{thm2.2}.
\qed

Next, we will derive Corollary \ref{cor2.3} from the estimate (\ref{2.4}). The first part of Corollary \ref{cor2.3} would
follow if we find an integer valued sequence $R_L,\, L=1,2,...$ such that $R_L\to\infty$, $\frac {R_L}{N_L}\to 0$ and
$N_L\phi(R_L-n_L)\to 0$ as $L\to\infty$. In order to do this we observe first that since $\phi(k)$ is summable and nonincreasing,
\[
[N/2]\phi(N)\leq\sum_{k=[N/2]}^N\phi(k)\leq\sum_{k=[N/2]}\phi(k)\to 0\quad\mbox{as}\quad N\to\infty
\]
which means that $\gam(N)=N\phi(N)\to 0$ as $N\to\infty$. Observe that if $\phi(k)=0$ for some $k\geq 1$ then by monotonicity
$\phi(n)=0$ for all $n\geq k$. In this case there is nothing to prove taking, say, $R_L=2n_L\to\infty$ as $L\to\infty$. Hence, we
can and will assume that $\phi(n)>0$ for all $n\geq 1$. For some $0<\ve<1$ set
\[
M_N=M_N^{(\ve)}=\min\{ n\geq 1:\, \frac n{\gam^\ve(n)}\geq N\}\to\infty\quad\mbox{as}\quad N\to\infty.
\]
Then
\[
\frac {M_N}{\gam^\ve(M_N)}\geq N\quad\mbox{and}\quad N\phi(M_N)=\frac N{M_N}\gam(M_N)\leq\gam^{1-\ve}(M_N)\to 0\,\,\mbox{as}\,\, N\to\infty.
\]
Since $\frac {M_N-1}{\gam^\ve(M_N-1)}< N$ then $\frac N{M_N-1}>\frac 1{\gam^\ve(M_N-1)}\to\infty$ as $N\to\infty$, and so $\frac {M_N}N\to 0$
as $N\to\infty$. Hence, taking $R_L=M_{N_L}+n_L$ we conclude the proof of the first part of Corollary \ref{cor2.3}.

In the second part of Corollary \ref{cor2.3} we set $V_{L}=A^\eta_{n_L}=[\eta_0,...,\eta_{n_L-1}]$ where $\eta$ is a nonperiodic sequence
and observe that $P(A^\eta_{n_L})\leq e^{-\upsilon n_L}$ by (\ref{2.3}). Hence, the conditions of the first part of Corollary \ref{cor2.3}
would hold true provided
\begin{equation}\label{6.14}
\pi(A^\eta_n)\to\infty\quad\mbox{as}\quad n\to\infty
\end{equation}
whenever $\eta$ is a nonperiodic sequence. To see this note that $\pi(A^\eta_n)$ is, clearly, nondecreasing in $n$, and so
$\lim_{n\to\infty}\pi(A^\eta_n)=r$ exists. If $r<\infty$ then there exists $n_0\geq 1$ such that $\pi(A^\eta_n)=r$ for all
$n\geq n_0$ which means that $\eta$ is periodic with the period $r$. Hence, $r=\infty$ since $\eta$ is not periodic completing
the proof of Corollary \ref{cor2.3}.
\qed

\section{Geometric distribution limits}\label{sec4}\setcounter{equation}{0}

It will be convenient to set $V^{(0)}=V\in\cF_{0,n-1}$, $V^{(1)}=W\in\cF_{0,m-1}$ and
\[
X^{(\al)}_{k,N}=\bbI_{V^{(\al)}}\circ T^{q_{N}(k)},\,\,\al=0,1
\]
 so that
 \[
 \tau=\tau_{V^{(1)}}=\min\{ k\geq 1:\, X^{(1)}_{k,N}=1\}\,\,\,\mbox{and}\,\,\,\Sig_N^{V^{(0)},V^{(1)}}=\sum_{k=1}^\tau X^{(0)}_{k,N}.
 \]
 Set also $S_L=\sum_{k=1}^LX_{k,N}^{(0)}$, so that $S_\tau=\Sig_N^{V^{(0)},V^{(1)}}$, and denote $\tau_N=\min(\tau,N)$. Let $\{ Y^{(\al)}_{k,N}:\,
 k\geq 1,\,\al=0,1\}$ be a sequence of independent Bernoulli random variables such that $Y^{(\al)}_{k,N}$ has the same distribution
 as $X^{(\al)}_{k,N}$. Since $P$ is $T$-invariant $E(X^{(\al)}_{k,N})=P\{ X^{(\al)}_{k,N}=1\}=E(Y^{(\al)}_{k,N})=P\{ Y^{(\al)}_{k,N}=1\}=P(V^{(\al)})$.
 Set
 \[
 S^*_L=\sum_{k=1}^LY^{(0)}_{k,N},\,\,\tau^*=\min\{ k\geq 1:\, Y^{(1)}_{k,N}=1\}\,\,\mbox{and}\,\,\tau^*_N=\min(\tau^*,N).
 \]
 We can and will assume that all above random variables are defined on the same (sufficiently large) probability space.
 By Lemma 3.1 from \cite{KR1} the sum $S^*_{\tau^*}$ has the geometric distribution with the parameter
 \begin{equation}\label{7.2}
 \vrho=\frac {P(V^{(1)})}{P(V^{(1)})+P(V^{(0)})(1-P(V^{(1)}))}>\rho
 \end{equation}
 where $\rho=P(V^{(1)})\big(P(V^{(1)})+P(V^{(0)})\big)^{-1}$.

 Next, we can write
\begin{equation}\label{7.3}
d_{TV}(\cL(S_\tau),\,\mbox{Geo}(\rho))\leq A_1+A_2+A_3+A_4
\end{equation}
where $A_1=d_{TV}(\cL(S_\tau),\,\cL(S_{\tau_N}))$,
 $A_2=d_{TV}(\cL(S_{\tau_N}),\,\cL(S^*_{\tau^*_N}))$,
 $A_3=d_{TV}(\cL(S^*_{\tau^*_N}),\,\cL(S^*_{\tau^*}))$ and
 $A_4=d_{TV}(\mbox{Geo}(\varrho),\,\mbox{Geo}(\rho))$.

 Introduce random vectors
 $\bfX_{N}^{(\al)}=\{ X_{k,n}^{(\al)},\, 1\leq k\leq N\},\,\al=0,1$, $\bfX_N=
 \{\bfX_{N}^{(0)},\,\bfX_{N}^{(1)}\}$, $\bfY_{N}^{(\al)}=\{ Y_{n,N}^{(\al)},\, 1\leq k\leq N\},
 \,\al=0,1$ and $\bfY_N=\{\bfY_{N}^{(0)},\,\bfY_{N}^{(1)}\}$.
  Observe that the event $\{ S_\tau\ne S_{\tau_N}\}$ can
 occur only if $\tau>N$. Also, we can write $\{\tau>N\}=\{ X_{n,N}^{(1)}=0\,\,
 \mbox{for all}\,\, k=1,...,N\}$ and $\{\tau^*>N\}=\{ Y_{n,0}^{(1)}=0\,\,
 \mbox{for all}\,\, k=1,...,N\}$ Hence,
 \begin{eqnarray}\label{7.4}
 &A_1\leq P\{\tau>N\}= P\{\tau^*>N\}+|P\{ X_{n,N}^{(1)}=0\,\,\mbox{for}\,\,
 n=1,...,N\}\\
 &-P\{ Y_{n,N}^{(1)}=0\,\,\mbox{for}\,\, n=0,1,...,N\}|\leq P\{\tau^*>N\}+d_{TV}(\cL(\bfX_{N}),\,\cL(\bfY_{N})).\nonumber
 \end{eqnarray}

 Since $Y^{(1)}_{k,N},\, k=0,1,...$ are i.i.d. random variables we obtain that
 \begin{equation}\label{7.5}
 P\{\tau^*>N\}=(1-P(V^{(1)}))^{N}.
 \end{equation}
 Also
 \begin{equation}\label{7.6}
 A_3\leq P\{\tau^*>N\}=(1-P(V^{(1)}))^{N}.
 \end{equation}
 The estimate of $A_4$ is also easy
 \begin{eqnarray}\label{7.7}
 &A_4\leq \sum_{k=0}^\infty |\vr(1-\vr)^k-\rho(1-\rho)^k|\leq 2\sum_{k=1}^\infty((1-\rho)^k-(1-\vr)^k)\\
 &=2(1-\rho)\rho^{-1}-2(1-\vr)\vr^{-1}=\frac {2(\vr-\rho)}{\rho\vr}=2P(V^{(1)}).\nonumber
 \end{eqnarray}

Next, we observe that by Theorem 3 in \cite{AGG},
\begin{equation}\label{7.9}
A_2\leq d_{TV}(\cL(\bfX_{N}),\,\cL(\bfY_{N}))\leq 2b_1+2b_2+b_3+2\sum_{1\leq k\leq N,\,\al=0,1}(p_{k,N}^{(\al)})^2
\end{equation}
where $p_{k,N}^{(\al)}=P\{ X_{k,N}^{(\al)}=1\}=P(V^{(\al)})$ while the definitions of $b_1,b_2$ and $b_3$ are similar to Section \ref{sec3}
taking into account the additional parameter $\al$. Namely, setting
\[
B^{R}_{k,N}=\{ (l,0),\, (l,1):\, 1\leq l\leq N,\,\del(k,l)< R\},\, p^{\al,\be}_{k,l,N}=E(X^{(\al)}_{k,N}X^{(\be)}_{l,N})
\]
and $I_N=\{(k,\al):\, 1\leq k\leq N,\,\al=0,1\}$ we have
\begin{equation}\label{7.11}
b_1=\sum_{(k,\al)\in I_N}\sum_{(l,\be)\in B^{R}_{k,N}}p^{(\al)}_{k,N}p^{(\be)}_{l,N},
\end{equation}
\begin{equation}\label{7.12}
b_2=\sum_{(k,\al)\in I_N}\sum_{(k,\al)\ne(l,\be)\in B^{R}_{k,N}}p^{(\al,\be)}_{k,l,N}\quad\mbox{and}
\end{equation}
\begin{equation}\label{7.13}
b_3=\sum_{(k,\al)\in I_N}s^{(\al)}_{k,N}\,\,\,\mbox{where}
\end{equation}
\[
s^{(\al)}_{k,N}=E\big\vert E\big(X^{(\al)}_{k,N}-p^{(\al)}_{k,N}|
\sig\{ X^{(\be)}_{l,N}:\, (l,\be)\in I_N\setminus B^{R}_{k,N}\}\big)\big\vert.
\]

Since $p^{(\al)}_{k,N}= P(V^{(\al)})$, it follows taking into account (\ref{3.3}) and (\ref{3.4}) that
\begin{equation}\label{7.15}
b_1\leq 6KRN((P(V^{(0)}))^2+(P(V^{(1)}))^2).
\end{equation}
In order to estimate $p^{\al,\be}_{k,l,N}$ (and, eventually, $b_2$) we observe that
\[
(\bbI_{V^{(0)}}\circ T^i)(\bbI_{V^{(1)}}\circ T^j)=0\,\,\mbox{if}\,\, |i-j|<\ka_{V^{(0)},V^{(1)}}.
\]
Hence, $p^{\al,\be}_{k,l,N}=0$ if $\del_N(k,l)<\ka_{V^{(0)},V^{(1)}}$. Now suppose that
$\del_N(k,l)=d\geq\ka_{V^{(0)},V^{(1)}}$. Then we have to deal with two alternatives from (\ref{6.5}).
 Assume, for instance, that the first inequality in (\ref{6.5}) holds true and let $r=q_{N}(k)-q_{N}(l)$.
Then $r\geq d\geq\ka_{V^{(0)},V^{(1)}}$. If $r\geq n$ then since $V^{(0)}\in\cF_{0,n-1}$ and $V^{(1)}\in\cF_{0,m-1}$,
we obtain by the definition of the $\phi$-mixing coefficient that
\begin{equation}\label{7.16}
p^{\al,\be}_{k,l,N}=P(V^{(\be)}\cap T^{-r}V^{(\al)})\leq(\phi(r-m+1)+P(V^{(\al)}))P(V^{(\be)}).
\end{equation}
Suppose that $\ka_{V^{(0)},V^{(1)}}\leq r<n$ and assume that $V^{(\be)}\cap T^{-r}V^{(\al)}\ne\emptyset$.
Let $s\geq n-r$, set $V_s=T^sV$ and observe that $T^{-s}V_s\supset V$. Then by the definition of the $\phi$-mixing
 coefficient,
\begin{eqnarray}\label{7.17}
&p^{\al,\be}_{k,l,N}=P(V^{(\be)}\cap T^{-r}V^{(\al)})\leq P(V^{(\be)}\cap T^{-(r+s)}V^{(\al)}_s)\\
&\leq (\phi(r+s-n+1)+P(T^sV^{(\al)}))P(V^{(\be)})\nonumber\\
&\leq (\phi([r/2]+1)+P(T^{n-[r/2]}V^{(\al)}))P(V^{(\be)})\nonumber
\end{eqnarray}
taking $s= n-[r/2]$. If the second inequality in (\ref{6.5}) holds true then we obtain (\ref{7.16}) if $r=q_{N}(l)-q_{N}(k)\geq n$,
  while if $\ka_{V^{(0)},V^{(1)}}\leq r<n$ then we arrive at (\ref{7.17}).
  and integers $k\geq 0$ and $r$,

  Now, it follows from (\ref{3.3}), (\ref{3.4}), (\ref{6.9}), (\ref{7.16}) and (\ref{7.17}) that
  \begin{eqnarray}\label{7.20}
&b_2=\sum_{k=1}^N\sum_{k\ne l\in B_{k,N}^R,\,\al,\be=0,1}p^{\al,\be}_{k,l,N}\leq 4KN(P(V^{(0)})+P(V^{(1)}))\\
&\times\sum_{r=\ka_{V^{(0)},V^{(1)}}}^R(\phi([r/2]+1)+P(T^{n-[r/2]}V^{(0)})+P(T^{m-[r/2]}V^{(1)})).\nonumber
\end{eqnarray}
Similarly to (\ref{6.13}) we obtain also that
\begin{equation}\label{7.21}
b_3=\sum_{1\leq k\leq N,\,\al=0,1}^Ns^{(\al)}_{k,N}\leq 6N\phi(R-n\vee m).
\end{equation}
These provide the estimate of $A_2$ by (\ref{7.9}), (\ref{7.11})--(\ref{7.15}), (\ref{7.20}) and (\ref{7.21}).
Finally, combining (\ref{7.3})--(\ref{7.15}), (\ref{7.20}) and (\ref{7.21}) we derive (\ref{2.6}) completing the
 proof of Theorem \ref{thm2.4}.       \qed

 Corollary \ref{cor2.5} follows from the estimate (\ref{2.6}) choosing $R=R_L$ as in Corollary \ref{cor2.3} and if $V_L=A_{n_L}^\xi$ and
 $W_L=A^\eta_{m_L}$ it remains
 only to verify the assertion that $\ka_{A_n^\xi,A_m^\eta}\to\infty$ as $n,m\to\infty$ provided that $\xi,\eta\in\Om_P$ are not periodic and not
 shifts of each other. Indeed, $\pi(A^\xi_n)$, $\pi(A^\eta_m)$ and $\pi(A^\xi_n,A^\eta_m)$ are nondecreasing in $n$ and $m$, and so does
 $\pi(A^\xi_n,A^\eta_m)$. Hence, the limit $r=\lim_{n,m\to\infty}\ka_{A_n^\xi,A_m^\eta}$ exists. If $r<\infty$ then, at least, one of the limits
 $r_1=\lim_{n\to\infty}\pi(A^\xi_n)$, $r_2=\lim_{m\to\infty}\pi(A^\eta_m)$ or $r_3=\lim_{n,m\to\infty}\pi(A^\xi_n,A^\eta_m)$ is finite. If $r_1<\infty$
 then $\xi$ is periodic with the period $r_1$, if $r_2<\infty$ then $\eta$ is periodic with the period $r_2$ and if $r_3<\infty$ then either $T^{r_3}\xi
 =\eta$ or $T^{r_3}\eta=\xi$. Finally, it follows from Lemma 3.2 from \cite{KY} that (\ref{2.13}) holds true for $P\times P$-almost all $(\xi,\eta)$,
   completing the proof.    \qed


\begin{thebibliography}{Bow75}

\itemsep=\smallskipamount



\bibitem{AV}
M. Abadi and N. Vergne, {\em Sharp errors for point-wise Poisson
approximations in mixing processes}, Nonlinearity 21 (2008), 2871--2885.



\bibitem{AGG}
R. Arratia, L. Goldstein and L. Gordon, {\em Two moments suffice for Poisson
approximations: the Chen--Stein method}, Ann. Probab. 17 (1989), 9--25.





\bibitem{Bo}
R. Bowen, {\em Equilibrium States and the Ergodic Theory of Anosov
Diffeomorphisms}, Lecture Notes in Math. 470, Springer--Verlag, Berlin, 1975.

\bibitem{Br} R.C. Bradley, {\em Introduction to Strong Mixing Conditions,}
Kendrick Press, Heber City, 2007.



\bibitem{DWY}
M. Demers, P.Wright and L.-S. Young, {\em Entropy, Lyapunov exponents and
escape rates in open systems}, Ergod. Th.\& Dynam. Sys. 30 (2012), 1270--1301.


\bibitem{De}
M. Denker, {\em Remarks on weak laws for fractal sets}, in: Fractal Geometry and
Stochastics (eds. C.Bandt, S.Graf and M.Z\" ahle), p.p.167--178, Birkh\" auser, Basel, 1995.


\bibitem{Ha} Ye. Hafouta, {\em  A functional CLT for nonconventional polynomial arrays },
arXiv: 1907.03303


\bibitem{He}
L. Heinrich, {\em Mixing properties and central limit theorem for a class of
non-identical piecewise monotonic $C^2$-transformations},
Mathematische Nachricht. 181 (1996), 185--214.

\bibitem{Hi}
M. Hirata, {\em Poisson law for Axiom A diffeomorphisms}, Ergod. Th. Dynam. Sys.
13 (1993), 533--556.

\bibitem{HK} Ye. Hafouta and Yu. Kifer, {\em Nonconventional Limit Theorems and
 Random Dynamics}, World Scientific, Singapore, 2018.


\bibitem{HW}
N.T.A. Haydn and K. Wasilewska, {\em Limiting distribution and error terms for the number of visits
to balls in non-uniform hyperbolic dynamical systems}, Discr. Cont. Dyn.Syst. 36 (2016), 2586--2611.


\bibitem{Ki1} Yu. Kifer, {\em Nonconventional Poisson limit theorems},
Israel J. Math. 195 (2013), 373--392.


\bibitem{Ki19} Yu. Kifer, {\em Limit theorems for numbers of multiple returns in
nonconventional arrays}, arXiv: 1910.01439.

\bibitem{KR1}
Yu. Kifer and A. Rapaport, {\em Poisson and compound Poisson approximations
in conventional and nonconventional setups}, Probab. Th. Relat. Fields 160
(2014), 797--831.

\bibitem{KR2} Yu. Kifer and A. Rapaport, {\em Geometric law for multiple returns
until a hazard}, Nonlinearity 32 (2019), 1525--1545.


\bibitem{KY} Yu. Kifer and F. Yang, {\em Geometric law for numbers of returns until
a hazard under $\phi$-mixing}, arXiv: 1812.09927.




\bibitem{Pi} B. Pitskel, {Poisson limit law for Markov chains}, Ergod. Th. Dynam. Sys.
11 (1991), 501--513.






\end{thebibliography}

\end{document}